\documentclass[11pt]{amsart}

\usepackage{amssymb,amsmath,amsthm}
\usepackage{enumerate}
\usepackage{booktabs}
\usepackage{tikz}
\usetikzlibrary{arrows.meta,calc}
\usepackage[colorlinks=true,linkcolor=blue,citecolor=blue,urlcolor=blue]{hyperref}

\newtheorem{theorem}{Theorem}[section]
\newtheorem{lemma}[theorem]{Lemma}
\newtheorem{proposition}[theorem]{Proposition}

\newtheorem{conjecture}[theorem]{Conjecture}
\newtheorem{mainthm}{Theorem}

\newtheorem{maincor}[mainthm]{Corollary}
\theoremstyle{definition}
\newtheorem{definition}[theorem]{Definition}
\newtheorem{example}[theorem]{Example}
\newtheorem{problem}{Problem}
\theoremstyle{remark}
\newtheorem{remark}[theorem]{Remark}

\newcommand{\R}{\mathbb{R}}

\newcommand{\I}{I}
\newcommand{\Fib}{\mathfrak{F}}
\newcommand{\eps}{\varepsilon}
\newcommand{\dist}{\operatorname{dist}}

\newcommand{\rk}{\operatorname{rk}}
\newcommand{\Ker}{\operatorname{Ker}}
\newcommand{\ICT}{\textup{ICT}}
\newcommand{\Z}{\mathbb{Z}}
\newcommand{\M}{\mathbf{M}}
\newcommand{\W}{\mathcal{W}}

\begin{document}

\title[$\omega$-limit sets of triangular maps over periodic orbits]{$\omega$-limit sets of triangular maps over periodic orbits of the base}

\author{Juan Luis Garc\'{\i}a Guirao}
\address{Departamento de Matem\'atica Aplicada y Estad\'{\i}stica,\newline
Universidad Polit\'ecnica de Cartagena, Hospital de Marina,
30203 Cartagena (Regi\'on de Murcia), Spain}
\email{juan.garcia@upct.es}

\subjclass[2020]{Primary 37E05, 37B20; Secondary 37B45, 54H20}
\keywords{Triangular map, skew product, $\omega$-limit set, internal chain transitivity, periodic orbit, Cantor--Bendixson rank}

\begin{abstract}
Let $F(x,y)=(f(x),g(x,y))$ be a continuous triangular map of the unit square and let $z=(x_0,y_0)$ be a point such that
$\omega_f(x_0)$ is a periodic orbit $\{p_0,\dots,p_{m-1}\}$ of $f$. Then
$\omega_F(z)=\bigcup_{i}\{p_i\}\times M_i$ for nonempty compact sets $M_i$, and we ask which $m$-tuples
$(M_0,\dots,M_{m-1})$ occur. For $m=1$ the answer is a theorem of Kolyada and Snoha. We show that a tuple occurs if and
only if there are continuous surjections $M_i\to M_{i+1}$, indices modulo $m$, whose composition around the cycle is chain
transitive. This is a cyclic version of classical facts due to Dowker and Friedlander, Sharkovsky, Hirsch, Smith and Zhao,
and Kolyada and Snoha. It yields two necessary conditions: no $M_i$ is one of the exceptional sets of Kolyada and Snoha,
and there are continuous surjections $M_0\to M_1\to\dots\to M_{m-1}\to M_0$. We prove that these conditions are also
sufficient when every $M_i$ has uncountably many connected components, and when every $M_i$ has countably many
components of which only finitely many are nondegenerate. In particular the sections of one $\omega$-limit set need not
be homeomorphic. Every tuple satisfying the two conditions falls into one of these cases or into a single remaining class,
in which the sections have countably many components and infinitely many of them are nondegenerate; that class is
studied in a companion paper.
\end{abstract}

\maketitle

\section{Introduction}\label{sec:intro}

Throughout, $\I=[0,1]$. A \emph{triangular map} (or \emph{skew product}) of the square $\I^2$ is a continuous map
\begin{equation}\label{eq:triangular}
F(x,y)=\bigl(f(x),\,g(x,y)\bigr),\qquad (x,y)\in\I^2 .
\end{equation}
The continuous map $f\colon\I\to\I$ is the \emph{base}, and for $x\in\I$ the map $g(x,\cdot)\colon\I\to\I$ is the
\emph{fibre map} over $x$. The map $F$ sends each \emph{fibre} $\I_x=\{x\}\times\I$ into $\I_{f(x)}$.
For a continuous map $\varphi\colon Y\to Y$ of a metric space and $y\in Y$, the \emph{$\omega$-limit set} $\omega_\varphi(y)$ is
the set of limit points of the sequence $(\varphi^n(y))_{n\ge0}$.
Since Kloeden's extension of Sharkovsky's theorem to triangular maps \cite{Kl} and Kolyada's paper \cite{Ko}, these maps
have served as a testing ground between interval dynamics and genuinely two-dimensional phenomena; see \cite{Sm}.

For interval maps the structure of $\omega$-limit sets is completely known.

\begin{theorem}[Agronsky, Bruckner, Ceder and Pearson \cite{ABCP}; Bruckner and Sm\'{\i}tal \cite{BS}]\label{thm:ABCP}
A nonempty closed set $C\subset\I$ is the $\omega$-limit set of some point under some continuous map $\I\to\I$ if and only if
$C$ is nowhere dense or $C$ is a finite union of nondegenerate closed intervals.
\end{theorem}

For triangular maps the corresponding problem is open. Jim\'enez L\'opez and Sm\'{\i}tal \cite{JS} characterized the
$\omega$-limit sets that are Peano continua and contain the orbit generating them (see \cite{JS} for these notions) and gave two similar non-Peano continua, only one of which is an
$\omega$-limit set; see also \cite{AC,Ba,JS2}. The case that is completely understood is that of $\omega$-limit sets
contained in one fibre.

\begin{definition}\label{def:exceptional}
A nonempty closed set $M\subset\R$ is \emph{KS-exceptional} if
$M=J_1\cup\dots\cup J_r\cup C$, where $r\ge1$, the $J_k$ are nondegenerate compact intervals, $C$ is nonempty and
countable, all these sets are pairwise disjoint, and $\dist(C,J_k)>0$ for at least one $k$; here
$\dist(A,B)=\inf\{|a-b|:a\in A,\ b\in B\}$.
We write $\Fib$ for the class of nonempty closed subsets of $\I$ that are not KS-exceptional.
\end{definition}

\begin{theorem}[Kolyada and Snoha {\cite[Theorem~1]{KS}}]\label{thm:KS}
Let $a\in\I$ and $M\subset\I$. There are a triangular map $F$ and a point $z\in\I^2$ with $\omega_F(z)=\{a\}\times M$
if and only if $M\in\Fib$.
\end{theorem}

Kolyada and Snoha also observed \cite[p.~116]{KS} that, for a finite set $B$, the product $B\times M$ is an $\omega$-limit
set of a triangular map if and only if $M\in\Fib$, and proved every closed $M\subset\I$ is the intersection of some $\omega$-limit
set with a fibre \cite[Theorem~2]{KS}. Kolyada and Robatian \cite{KR} studied the analogous question for triangular maps
induced on the hyperspace of subintervals; see also \cite{MaR}.

\medskip
\noindent\textbf{The problem.}
A set $P=\{p_0,\dots,p_{m-1}\}\subset\I$ of $m$ distinct points is a \emph{periodic orbit} of $f$ if
$f(p_i)=p_{i+1}$, indices taken in $\Z_m=\Z/m\Z$. Let $z=(x_0,y_0)$ and suppose that $\omega_f(x_0)$ is finite. Then
$\omega_f(x_0)$ is a periodic orbit $P$ and $f^{km+i}(x_0)\to p_{i+\nu}$ as $k\to\infty$ for some fixed $\nu\in\Z_m$ (see Lemma~\ref{lem:phase}). Since the first
coordinate of $F^n(z)$ is $f^n(x_0)$, every limit point of the orbit of $z$ lies in a fibre over $P$, and each fibre
$\I_{p_j}$ contains limit points of the subsequence $(F^{km+j-\nu}(z))_k$; hence
\begin{equation}\label{eq:tuple}
\omega_F(z)=\bigcup_{i\in\Z_m}\{p_i\}\times M_i,\qquad M_i\subset\I\ \text{nonempty and closed}.
\end{equation}
This paper studies which tuples $(M_0,\dots,M_{m-1})$ occur.

\begin{definition}\label{def:realizable}
An $m$-tuple $\M=(M_0,\dots,M_{m-1})$ of nonempty closed subsets of $\I$ is \emph{realizable} if there are a triangular map
$F=(f,g)$, a periodic orbit $P=\{p_0,\dots,p_{m-1}\}$ of $f$ with $f(p_i)=p_{i+1}$, and a point $z=(x_0,y_0)$ such that
$\omega_f(x_0)=P$ and \eqref{eq:tuple} holds.
\end{definition}

For $m=1$, Theorem~\ref{thm:KS} says that $(M_0)$ is realizable if and only if $M_0\in\Fib$, and the result on products
covers constant tuples $M_0=\dots=M_{m-1}$.

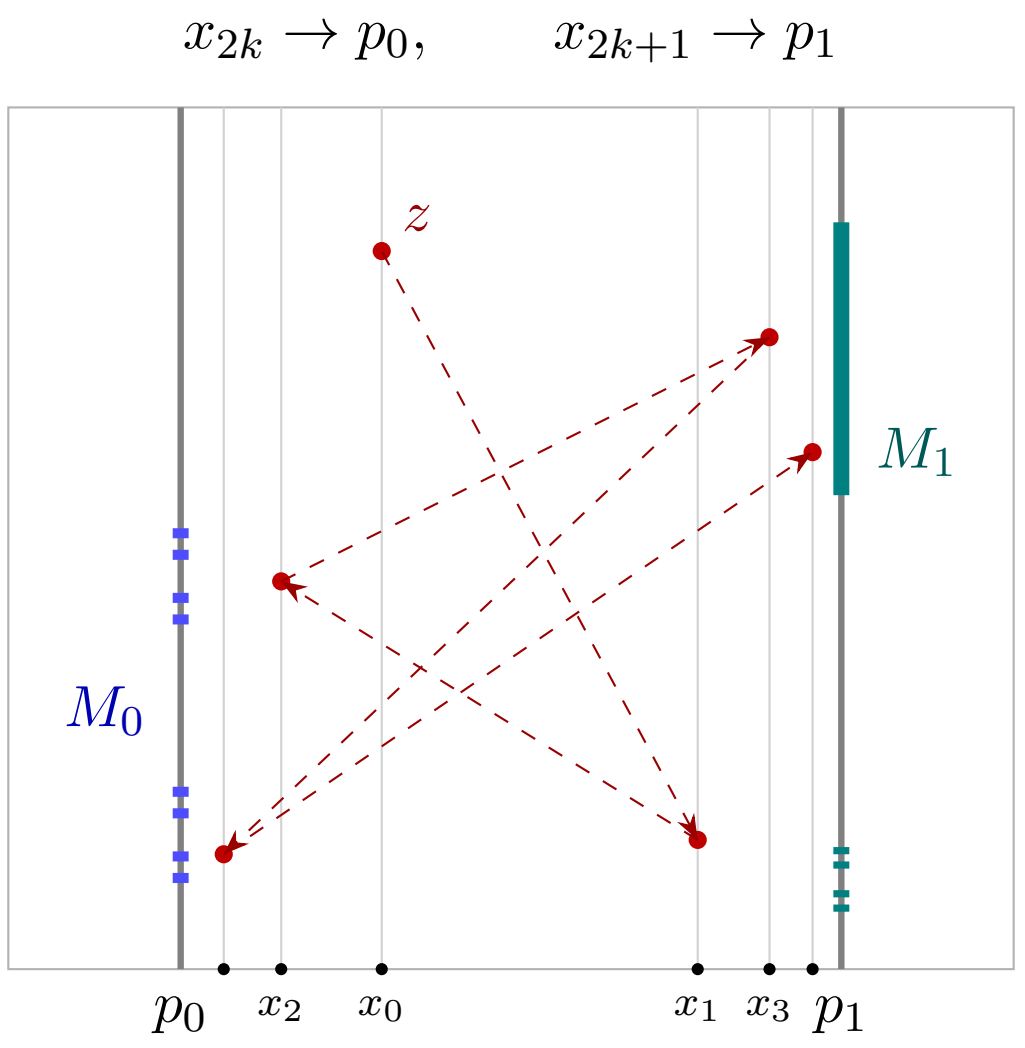
\begin{figure}[t]
\centering
\begin{tikzpicture}[scale=0.95,>=Stealth]
  \draw[gray!60] (0,0) rectangle (7,6);
  \draw[very thick,gray] (1.2,0)--(1.2,6);
  \draw[very thick,gray] (5.8,0)--(5.8,6);
  \node[below] at (1.2,-0.05) {$p_0$};
  \node[below] at (5.8,-0.05) {$p_1$};
  \foreach \y in {0.6,0.75,1.05,1.2,2.4,2.55,2.85,3.0}{\draw[line width=3pt,blue!70] (1.2,\y)--(1.2,\y+0.07);}
  \node[blue!70!black,left] at (1.1,1.8) {$M_0$};
  \foreach \y in {0.4,0.5,0.7,0.8}{\draw[line width=3pt,teal] (5.8,\y)--(5.8,\y+0.05);}
  \draw[line width=3pt,teal] (5.8,3.3)--(5.8,5.2);
  \node[teal!70!black,right] at (5.9,3.6) {$M_1$};
  \foreach \x in {2.6,4.8,1.9,5.3,1.5,5.6}{
     \draw[gray!35] (\x,0)--(\x,6);
     \fill (\x,0) circle (1.2pt);}
  \node[below] at (2.6,-0.05) {$\scriptstyle x_0$};
  \node[below] at (4.8,-0.05) {$\scriptstyle x_1$};
  \node[below] at (1.9,-0.05) {$\scriptstyle x_2$};
  \node[below] at (5.3,-0.05) {$\scriptstyle x_3$};
  \coordinate (z0) at (2.6,5.0);
  \coordinate (z1) at (4.8,0.9);
  \coordinate (z2) at (1.9,2.7);
  \coordinate (z3) at (5.3,4.4);
  \coordinate (z4) at (1.5,0.8);
  \coordinate (z5) at (5.6,3.6);
  \foreach \p in {z0,z1,z2,z3,z4,z5}{\fill[red!75!black] (\p) circle (1.8pt);}
  \draw[->,red!60!black,dashed] (z0)--(z1);
  \draw[->,red!60!black,dashed] (z1)--(z2);
  \draw[->,red!60!black,dashed] (z2)--(z3);
  \draw[->,red!60!black,dashed] (z3)--(z4);
  \draw[->,red!60!black,dashed] (z4)--(z5);
  \node[above right,red!60!black] at (z0) {$z$};
  \node at (3.5,6.45) {$x_{2k}\to p_0,\qquad x_{2k+1}\to p_1$};
\end{tikzpicture}
\caption{An $\omega$-limit set over a periodic orbit $\{p_0,p_1\}$ of the base. The orbit of $z$ alternates between
neighbourhoods of the fibres $\I_{p_0}$ and $\I_{p_1}$, and $\omega_F(z)=\{p_0\}\times M_0\cup\{p_1\}\times M_1$. The
two sections may be topologically different; in the picture $M_0$ is a Cantor set and $M_1$ is a Cantor set together with
an interval, as in Example~\ref{ex:cantor}.}
\label{fig:setting}
\end{figure}

\medskip
\noindent\textbf{Chains.}
Let $\psi\colon X\to X$ be a continuous map of a compact metric space $(X,d)$ and let $\eps>0$. An \emph{$\eps$-chain} from
$a$ to $b$ is a finite sequence $a=x_0,x_1,\dots,x_n=b$ with $n\ge1$ and $d(\psi(x_{k-1}),x_k)<\eps$ for $1\le k\le n$. The
map $\psi$ is \emph{chain transitive} if for all $a,b\in X$ and $\eps>0$ there is an $\eps$-chain from $a$ to $b$. A closed set
$K\subset X$ is \emph{internally chain transitive} (\ICT) for $\psi$ if such chains exist for all $a,b\in K$ and $\eps>0$
with all $x_k\in K$.

Given an $m$-tuple $\M$ of nonempty compact subsets of $\I$, let $X_\M=\bigcup_{i\in\Z_m}\{i\}\times M_i$ be their disjoint
union. On $X_\M$ and on $\Z_m\times\I$ we use the metric that agrees with the Euclidean one inside each piece
$\{i\}\times\I$ and equals $2$ between different pieces; by Lemma~\ref{lem:metric} below, chain transitivity does not depend on
this choice. A
continuous map $\Psi\colon X_\M\to X_\M$ is \emph{cyclic} if $\Psi(\{i\}\times M_i)\subset\{i+1\}\times M_{i+1}$ for every $i$.
Given continuous maps $\varphi_0,\dots,\varphi_{m-1}\colon\I\to\I$, the \emph{induced cyclic map} on $\Z_m\times\I$ is
$\Phi(i,y)=(i+1,\varphi_i(y))$.

\medskip
\noindent\textbf{Main results.}

\begin{mainthm}[Realization over a cycle]\label{thm:A}
For an $m$-tuple $\M$ of nonempty closed subsets of $\I$ the following are equivalent:
\begin{enumerate}[\upshape(a)]
\item $\M$ is realizable;
\item $X_\M$ admits a chain transitive cyclic map;
\item there are continuous surjections $\psi_i\colon M_i\to M_{i+1}$, $i\in\Z_m$, such that the return map
$R=\psi_{m-1}\circ\dots\circ\psi_0\colon M_0\to M_0$ is chain transitive.
\end{enumerate}
Moreover, the fibre maps over $P$ can be prescribed. Let $f\colon\I\to\I$ be continuous with a periodic orbit
$P=\{p_0,\dots,p_{m-1}\}$, $f(p_i)=p_{i+1}$, let $x_0\in\I$ satisfy $f^{km+i}(x_0)\to p_i$ as $k\to\infty$ for every $i$ and
$f^n(x_0)\notin P$ for all $n\ge0$, and let $\varphi_0,\dots,\varphi_{m-1}\colon\I\to\I$ be continuous. Then the tuples $\M$ for
which there are a continuous $g\colon\I^2\to\I$ with $g(p_i,\cdot)=\varphi_i$ and a point $y_0$ with
$\omega_{(f,g)}(x_0,y_0)=\bigcup_i\{p_i\}\times M_i$ are exactly those for which $\bigcup_i\{i\}\times M_i$ is a closed
\ICT\ set of the induced cyclic map $\Phi$.
\end{mainthm}

Theorem~\ref{thm:A} is a cyclic version of known facts, and we make the attributions precise. An \emph{asymptotic
pseudo-orbit} of $\psi$ is a sequence $(u_n)$ in $X$ with $d(u_{n+1},\psi(u_n))\to0$. That the set of accumulation points of an
asymptotic pseudo-orbit is \ICT, and conversely, is due to Hirsch, Smith and Zhao \cite{HSZ}; see also Barwell, Good,
Oprocha and Raines \cite{BGOR}, who proved that \ICT\ coincides with the weak incompressibility (see \cite{BGOR} for the
definition) of Sharkovsky \cite{Sh} and Dowker and Friedlander \cite{DF}, and Meddaugh and Raines \cite{MR} for the relation
between \ICT\ sets and $\omega$-limit sets of interval maps. Dowker and Friedlander and Sharkovsky showed that weakly incompressible sets are
$\omega$-limit sets of points of systems containing them; see \cite[Section~3]{KR}. For $m=1$, the implication from (b) to (a),
including the use of an extension theorem, is Lemma~2 of \cite{KS}, and the combination of Theorem~\ref{thm:A} with
Theorem~\ref{thm:KS} recovers the fact, implicit in the proof of \cite[Theorem~1]{KS}, that a nonempty compact set $K\subset\R$
admits a chain transitive continuous self-map if and only if it is not KS-exceptional (Remark~\ref{rem:KS}). What is new is
the cyclic setting and its consequences for tuples of different sets.

\begin{maincor}[Necessary conditions]\label{cor:B}
If $\M=(M_0,\dots,M_{m-1})$ is realizable, then:
\begin{enumerate}[\upshape(i)]
\item $M_i\in\Fib$ for every $i$;
\item there are continuous surjections $M_0\to M_1\to\dots\to M_{m-1}\to M_0$.
\end{enumerate}
\end{maincor}

\begin{conjecture}\label{conj}
An $m$-tuple of nonempty closed subsets of $\I$ is realizable if and only if it satisfies conditions \textup{(i)} and
\textup{(ii)} of Corollary~\ref{cor:B}.
\end{conjecture}

We prove the conjecture outside one class of tuples. A \emph{component} of a set is a connected component; a component of a
compact subset of $\R$ is a point or a compact interval, and it is \emph{nondegenerate} if it is an interval.

\begin{mainthm}\label{thm:C}
Let $\M$ be an $m$-tuple of nonempty closed subsets of $\I$, each of which has uncountably many components. Then $\M$ is
realizable.
\end{mainthm}

No compatibility between consecutive sections is needed in Theorem~\ref{thm:C}; for instance a Cantor set may be followed
by a Cantor set together with an interval (Example~\ref{ex:cantor}).

\begin{mainthm}\label{thm:D}
Let $\M$ be an $m$-tuple of nonempty closed subsets of $\I$, each of which has countably many components, only finitely many
of them nondegenerate. Then $\M$ is realizable if and only if it satisfies conditions \textup{(i)} and \textup{(ii)} of
Corollary~\ref{cor:B}.
\end{mainthm}

The sections of a realizable tuple as in Theorem~\ref{thm:D} need not be homeomorphic (Example~\ref{ex:nonhomeo}), unlike
in the following special case.

\begin{maincor}\label{cor:E}
Let $\M$ be an $m$-tuple of nonempty closed subsets of $\I$, each with finitely many components. Then $\M$ is realizable if
and only if the sets $M_0,\dots,M_{m-1}$ are pairwise homeomorphic and each of them is finite or a finite union of
nondegenerate closed intervals.
\end{maincor}

Let $\W$ be the class of closed subsets of $\I$ that have countably many components, infinitely many of which are
nondegenerate. Every set in $\W$ belongs to $\Fib$, since a KS-exceptional set has only finitely many nondegenerate
components.

\begin{mainthm}\label{thm:F}
Let $\M$ satisfy condition \textup{(ii)} of Corollary~\ref{cor:B}. Then either every $M_i$ has uncountably many components,
or every $M_i$ has countably many components and they all have the same number of nondegenerate components, finite or
infinite. Consequently Conjecture~\ref{conj} holds for every tuple except possibly when every $M_i$ belongs to $\W$, in
which case it asserts that condition \textup{(ii)} alone implies realizability.
\end{mainthm}

The class $\W$ is studied in the companion paper \cite{GuII}. There we show that the realizable tuples in $\W$ need not
be homeomorphic either, that two Cantor--Bendixson ranks attached to the sections are invariant under (ii), and that
Conjecture~\ref{conj} holds for large subclasses of $\W$; for instance it holds when every degenerate component of every
section is a limit of nondegenerate components.

The paper is organized as follows. Section~\ref{sec:prelim} collects facts on chains, on countable compacta and on
components. Section~\ref{sec:A} proves Theorem~\ref{thm:A} and Corollary~\ref{cor:B}. Theorems~\ref{thm:C}
and~\ref{thm:D}, Corollary~\ref{cor:E} and Theorem~\ref{thm:F} are proved in Sections~\ref{sec:C} and~\ref{sec:D}.
Section~\ref{sec:examples} gives examples, and Section~\ref{sec:problems} contains remarks and open problems.

\section{Preliminaries}\label{sec:prelim}

\subsection{Chains}
Throughout this subsection $(X,d)$ is a compact metric space and $\psi\colon X\to X$ is continuous. For a sequence
$(u_n)_{n\ge0}$ in $X$ we write $\omega((u_n))$ for its set of accumulation points. The sequence is an \emph{asymptotic
pseudo-orbit} of $\psi$ if $\delta_n:=d(u_{n+1},\psi(u_n))\to0$. On a disjoint union $\bigcup_{i\in\Z_m}\{i\}\times Y_i$ of
compact subsets of $\I$ we use, as in Section~\ref{sec:intro}, the metric that agrees with the Euclidean one inside each piece and equals $2$ between
different pieces. For $p\in\R$ and $r>0$ we write $B(p,r)=(p-r,p+r)\cap\I$. The \emph{Hausdorff distance} between
nonempty compact sets $A,B\subset\R$ is the infimum of the $\eps>0$ such that each of them lies in the $\eps$-neighbourhood
of the other; a sequence of nonempty compact sets $Y_n$ \emph{converges} to a point $a$, written $Y_n\to a$, if every
neighbourhood of $a$ contains $Y_n$ for all but finitely many $n$, that is, if $Y_n\to\{a\}$ in the Hausdorff metric.

\begin{lemma}\label{lem:invariant}
If $K\subset X$ is closed and \ICT\ for $\psi$, then $\psi(K)=K$.
\end{lemma}

\begin{proof}
Let $x\in K$. For every $\eps>0$ an $\eps$-chain in $K$ starting at $x$ has a second point $x_1\in K$ with
$d(\psi(x),x_1)<\eps$; hence $\psi(x)\in K$. Conversely, let $b\in K$. For every $\eps>0$ an $\eps$-chain in $K$ ending at $b$
has a penultimate point $x\in K$ with $d(\psi(x),b)<\eps$, so $b\in\overline{\psi(K)}=\psi(K)$.
\end{proof}

The next two lemmas are due to Hirsch, Smith and Zhao \cite{HSZ}; see also \cite{BGOR}. We include short proofs because we
use the specific form of the second one.

\begin{lemma}[\cite{HSZ}]\label{lem:pseudo-ICT}
If $(u_n)$ is an asymptotic pseudo-orbit of $\psi$, then $K=\omega((u_n))$ is nonempty, closed and \ICT\ for $\psi$.
\end{lemma}

\begin{proof}
Nonemptiness and closedness follow from compactness. Let $a,b\in K$ and $\eps>0$. Choose $0<\eta<\eps/3$ such that
$d(u,v)<\eta$ implies $d(\psi(u),\psi(v))<\eps/3$. By compactness $\dist(u_n,K)\to0$. Choose $N$ with $\delta_n<\eps/3$ and
$\dist(u_n,K)<\eta$ for $n\ge N$, and pick $\ell>n_0\ge N$ with $d(u_{n_0},a)<\eta$ and $d(u_\ell,b)<\eta$. For $n_0<j<\ell$
pick $k_j\in K$ with $d(k_j,u_j)<\eta$, and set $k_{n_0}=a$, $k_\ell=b$. For $n_0\le j<\ell$,
\begin{align*}
d(\psi(k_j),k_{j+1})&\le d(\psi(k_j),\psi(u_j))+d(\psi(u_j),u_{j+1})+d(u_{j+1},k_{j+1})\\
&<\tfrac{\eps}{3}+\tfrac{\eps}{3}+\eta<\eps ,
\end{align*}
so $k_{n_0},\dots,k_\ell$ is an $\eps$-chain in $K$ from $a$ to $b$.
\end{proof}

\begin{lemma}[\cite{HSZ}]\label{lem:ICT-pseudo}
Let $K\subset X$ be closed and \ICT\ for $\psi$, and let $a_1\in K$. There is an asymptotic pseudo-orbit $(u_n)$ of $\psi$
with $u_0=a_1$, $u_n\in K$ for all $n$, $\omega((u_n))=K$, and $d(\psi(u_n),u_{n+1})<1$ for all $n$.
\end{lemma}

\begin{figure}[t]
\centering
\begin{tikzpicture}[scale=1,>=Stealth]
  \draw[thick,blue!70,fill=blue!5] (0,0) ellipse (4 and 1.5);
  \node[blue!70!black] at (3.3,1.45) {$K$};
  \fill (-2.8,0.1) circle (2pt); \node[left] at (-2.8,0.1) {$a_1$};
  \fill (-0.6,0.9) circle (2pt); \node[above] at (-0.6,0.9) {$a_2$};
  \fill (1.7,-0.4) circle (2pt); \node[right] at (1.7,-0.4) {$a_3$};
  \fill (0.2,-1.0) circle (2pt); \node[below] at (0.2,-1.0) {$a_4$};
  \foreach \a/\b in {{(-2.8,0.1)}/{(-2.45,0.5)},{(-2.45,0.5)}/{(-2.0,0.8)},{(-2.0,0.8)}/{(-1.5,0.95)},{(-1.5,0.95)}/{(-1.0,0.95)},
     {(-1.0,0.95)}/{(-0.6,0.9)},{(-0.6,0.9)}/{(-0.1,0.75)},{(-0.1,0.75)}/{(0.45,0.55)},{(0.45,0.55)}/{(0.95,0.3)},
     {(0.95,0.3)}/{(1.4,0.0)},{(1.4,0.0)}/{(1.7,-0.4)},{(1.7,-0.4)}/{(1.35,-0.8)},{(1.35,-0.8)}/{(0.8,-1.0)},
     {(0.8,-1.0)}/{(0.2,-1.0)},{(0.2,-1.0)}/{(-0.5,-0.95)},{(-0.5,-0.95)}/{(-1.2,-0.8)},{(-1.2,-0.8)}/{(-1.85,-0.55)},
     {(-1.85,-0.55)}/{(-2.4,-0.25)},{(-2.4,-0.25)}/{(-2.8,0.1)}}{
     \draw[->,red!70!black,shorten >=1.5pt,shorten <=1.5pt] \a--\b;}
  \foreach \p in {(-2.45,0.5),(-2.0,0.8),(-1.5,0.95),(-1.0,0.95),(-0.1,0.75),(0.45,0.55),(0.95,0.3),(1.4,0.0),
                  (1.35,-0.8),(0.8,-1.0),(-0.5,-0.95),(-1.2,-0.8),(-1.85,-0.55),(-2.4,-0.25)}{
     \fill[red!70!black] \p circle (1pt);}
  \node[red!70!black,align=center] at (0,-2.15)
     {stage $k$: $\frac1{k+1}$-chains in $K$ through $a_1,\dots,a_{k+1}$ and back to $a_1$};
\end{tikzpicture}
\caption{Proof of Lemma~\ref{lem:ICT-pseudo}. Concatenating the stages gives a sequence in $K$ whose one-step errors tend
to $0$ and which visits every $a_j$ infinitely often.}
\label{fig:stages}
\end{figure}
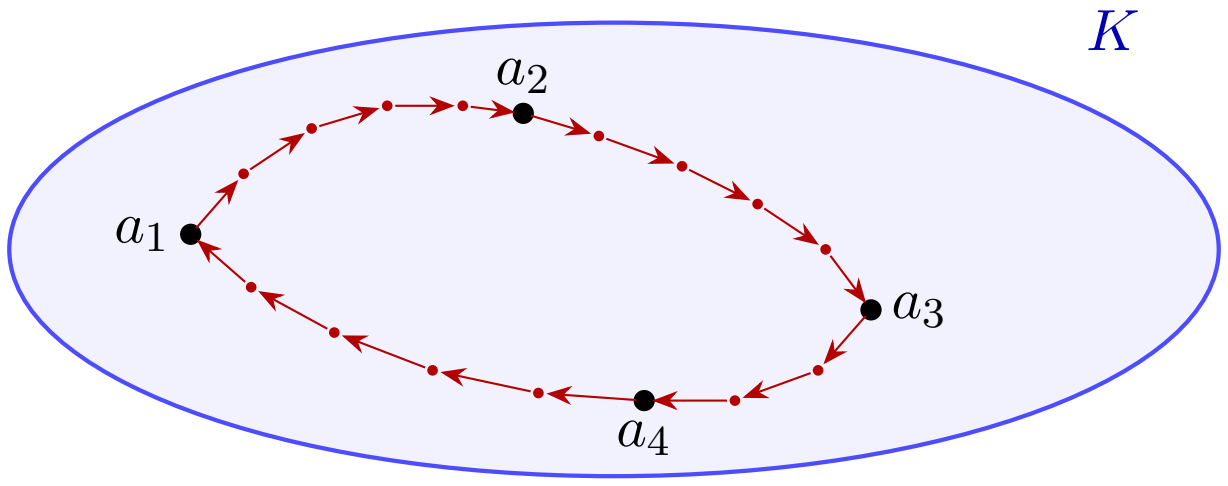

\begin{proof}
Extend $a_1$ to a sequence $(a_j)_{j\ge1}$ dense in $K$ (if $K$ is finite, let it enumerate $K$ cyclically). For $k\ge1$ let
the \emph{stage} $S_k$ be the concatenation of $\frac1{k+1}$-chains in $K$ from $a_1$ to $a_2$, from $a_2$ to $a_3$, \dots,
from $a_{k+1}$ back to $a_1$, writing each common endpoint once. Let $(u_n)$ be the concatenation $S_1S_2S_3\cdots$, again
identifying the final $a_1$ of each stage with the initial $a_1$ of the next (Figure~\ref{fig:stages}). Inside stage $k$
every step has error $<\frac1{k+1}\le\frac12$, so $\delta_n\to0$ and $d(\psi(u_n),u_{n+1})<1$. Every $a_j$ occurs in all
stages $k\ge j$, hence $K\subset\omega((u_n))$; conversely $u_n\in K$ and $K$ is closed.
\end{proof}

\begin{lemma}\label{lem:metric}
Chain transitivity of a continuous self-map of a compact metrizable space does not depend on the compatible metric.
Consequently, if $h\colon K\to K'$ is a homeomorphism and $\psi$ is chain transitive on $K$, then $h\circ\psi\circ h^{-1}$ is chain
transitive on $K'$.
\end{lemma}

\begin{proof}
Two compatible metrics on a compact space are uniformly equivalent, so for every $\eps>0$ there is $\eta>0$ such that every
$\eta$-chain for the first metric is an $\eps$-chain for the second. The second claim follows by transporting the metric.
\end{proof}

\begin{lemma}\label{lem:blocks}
For every $\eps>0$ and $m\ge1$ there is $\eta>0$ such that every $\eta$-chain $x_0,\dots,x_m$ of $\psi$ satisfies
$d(\psi^m(x_0),x_m)<\eps$.
\end{lemma}

\begin{proof}
Put $\eps_m=\eps$ and, going down, choose $0<\eps_k\le\eps_{k+1}/2$ such that $d(u,v)<\eps_k$ implies
$d(\psi(u),\psi(v))<\eps_{k+1}/2$. Let $\eta=\eps_1$. By induction $d(\psi^k(x_0),x_k)<\eps_k$: this holds for $k=1$, and
\[
d(\psi^{k+1}(x_0),x_{k+1})\le d\bigl(\psi(\psi^{k}(x_0)),\psi(x_k)\bigr)+d(\psi(x_k),x_{k+1})<\tfrac{\eps_{k+1}}2+\eta\le\eps_{k+1}. \qedhere
\]
\end{proof}

\begin{lemma}\label{lem:interval-CT}
Let $J$ be a nondegenerate compact interval. The identity of $J$ and the affine reflection of $J$ that exchanges its
endpoints are chain transitive.
\end{lemma}

\begin{proof}
Let $x,y\in J$ and $\eps>0$. For the identity, points $x=x_0,x_1,\dots,x_n=y$ on the segment from $x$ to $y$ with
consecutive distances $<\eps$ (and $n\ge1$) form an $\eps$-chain. For the reflection $\rho$, choose $n\ge1$ and $s$ with
$|s|<\eps$ and $x+ns=y$, and put $x_{2k}=x+ks$ and $x_{2k+1}=\rho(x_{2k})$ for $0\le k<n$, and $x_{2n}=y$. Then
$x_{2k+1}=\rho(x_{2k})$ exactly and $|\rho(x_{2k+1})-x_{2k+2}|=|x_{2k}-x_{2k+2}|=|s|<\eps$.
\end{proof}

\subsection{Countable compacta}
For a topological space $Y$ let $Y'$ be its set of accumulation points, and define the \emph{Cantor--Bendixson derivatives}
$Y^{(0)}=Y$, $Y^{(\beta+1)}=(Y^{(\beta)})'$ and $Y^{(\lambda)}=\bigcap_{\beta<\lambda}Y^{(\beta)}$ for limit ordinals $\lambda$.
If $Y$ is a nonempty countable compact metric space, there is a unique countable ordinal $\beta$ such that $Y^{(\beta)}$ is
finite and nonempty \cite[Section~6.C]{Ke}; we call $\beta$ the \emph{rank} of $Y$, write $\rk Y=\beta$, and call the finite set
$\Ker Y=Y^{(\beta)}$ the \emph{kernel} of $Y$. The \emph{rank} $\rk_Y(y)$ of a point $y\in Y$ is the largest $\gamma$ with
$y\in Y^{(\gamma)}$. If $W\subset Y$ is clopen, then $W^{(\gamma)}=W\cap Y^{(\gamma)}$ for every $\gamma$. A countable compact
metric space is zero-dimensional, so it has a basis of clopen sets.

\begin{theorem}[Mazurkiewicz and Sierpi\'nski \cite{MS}]\label{thm:MS}
Two nonempty countable compact metric spaces are homeomorphic if and only if they have the same rank and their kernels have
the same number of points.
\end{theorem}

\begin{lemma}\label{lem:CB}
Let $\pi\colon Y\to Z$ be a continuous surjection between compact metric spaces. Then $Z^{(\beta)}\subset\pi(Y^{(\beta)})$ for
every ordinal $\beta$. If $Y$ and $Z$ are countable, then $\rk Z\le\rk Y$, and $|\Ker Z|\le|\Ker Y|$ when the ranks are equal.
\end{lemma}

\begin{proof}
We use induction on $\beta$; the case $\beta=0$ is the surjectivity of $\pi$. Let $z\in Z^{(\beta+1)}$ and let $(z_n)$ be a
sequence of distinct points of $Z^{(\beta)}$ converging to $z$. Choose $y_n\in Y^{(\beta)}$ with $\pi(y_n)=z_n$; the $y_n$
are distinct, and a subsequence converges to some $y$, which belongs to $Y^{(\beta+1)}$ and satisfies $\pi(y)=z$. For a limit
ordinal $\lambda$ and $z\in Z^{(\lambda)}$, the sets $\pi^{-1}(z)\cap Y^{(\beta)}$, $\beta<\lambda$, are nonempty, compact and
decreasing, so their intersection is nonempty. The last assertion follows by taking $\beta=\rk Y+1$ and $\beta=\rk Y$.
\end{proof}

\begin{theorem}[{\cite[Lemma~8]{KS}}, see also {\cite[Lemma~6]{BS}}]\label{thm:KS8}
Let $A,B\subset\I$ be nonempty countable compact sets with $\rk A\ge\rk B$ and $\Ker B=\{b\}$. Then there is a continuous
surjection $\pi\colon A\to B$ with $\pi(\Ker A)=\{b\}$.
\end{theorem}

In \cite{KS} the rank is called the \emph{depth} and the kernel is denoted in the same way.

We say that a nonempty countable compact set $A$ \emph{dominates} a nonempty countable compact set $B$ if $\rk A>\rk B$, or
$\rk A=\rk B$ and $|\Ker A|\ge|\Ker B|$.

\begin{lemma}\label{lem:dominate}
If $A,B\subset\I$ are nonempty countable compact sets and $A$ dominates $B$, then there is a continuous surjection $A\to B$.
\end{lemma}

\begin{proof}
Let $\beta=\rk B$ and $\Ker B=\{b_1,\dots,b_c\}$. Since $B$ is zero-dimensional, we can write $B=B_1\cup\dots\cup B_c$ with
pairwise disjoint clopen sets $B_k\ni b_k$; then $\rk B_k=\beta$ and $\Ker B_k=\{b_k\}$. If $\rk A>\beta$, then
$A^{(\beta)}\supset A^{(\beta+1)}\neq\emptyset$ is infinite; if $\rk A=\beta$, then $|A^{(\beta)}|\ge c$. In both cases we can
pick distinct $a_1,\dots,a_c\in A^{(\beta)}$ and write $A=A_1\cup\dots\cup A_c$ with pairwise disjoint clopen sets
$A_k\ni a_k$; then $\rk A_k\ge\beta$. Theorem~\ref{thm:KS8} gives continuous surjections $A_k\to B_k$, and their union is a
continuous surjection $A\to B$.
\end{proof}

\subsection{Sets with uncountably many components}
A \emph{Cantor-like set} is a nonempty, nowhere dense, perfect subset of $\R$. If $P$ is a Cantor-like set, an
\emph{interval contiguous to $P$} is a bounded component of $\R\setminus P$. A point $a\in P$ is a \emph{bilateral
condensation point} of $P$ if $(a-\delta,a)\cap P$ and $(a,a+\delta)\cap P$ are uncountable for every $\delta>0$; every
point of $P$ other than $\min P$, $\max P$ and the countably many endpoints of contiguous intervals has this property.

\begin{theorem}[{\cite[Lemma~4]{KS}}]\label{thm:KS4}
Let $B\subset\I$ be a nonempty compact set and let $A\subset\I$ be a compact set containing a Cantor-like set $P$ such that
no interval contiguous to $P$ is a subset of $A$. Then there is a continuous surjection $A\to B$.
\end{theorem}

As in Section~\ref{sec:prelim}, $Y_n\to a$ means that every neighbourhood of $a$ contains $Y_n$ for all but finitely many
$n$.

\begin{theorem}[{\cite[Lemma~6 and the definition of homoclinic sets in Section~4]{KS}}]\label{thm:KS6}
Let $M\subset\I$ be a compact set containing a Cantor-like set $P$ such that no interval contiguous to $P$ is a subset of $M$,
and let $a$ be a bilateral condensation point of $P$. Then there are pairwise disjoint compact intervals $K_n$, $n\ge0$, not
containing $a$, such that $M\setminus\bigcup_nK_n=\{a\}$, the sets $M^n=M\cap K_n$ are nonempty, $M^n\to a$, and each
$M^n$ contains a Cantor-like set $P_n$ such that no interval contiguous to $P_n$ is a subset of $M^n$.
\end{theorem}

We also use the theorem of Alexandrov and Hausdorff that every uncountable Borel subset of $\R$ contains a Cantor-like set
\cite[Theorem~13.6]{Ke}.

\begin{lemma}[{cf.\ \cite[Lemmas~6 and~7]{KS}}]\label{lem:uncountable}
Let $M\subset\I$ be a compact set with uncountably many components. There are $a\in M$ and pairwise disjoint nonempty
clopen subsets $M^n$, $n\ge0$, of $M$ such that $M=\{a\}\cup\bigcup_{n\ge0}M^n$, $M^n\to a$, and every $M^n$ admits a
continuous surjection onto every nonempty compact subset of $\I$.
\end{lemma}

\begin{proof}
We repeat the argument of \cite[Lemma~7]{KS}. Only countably many components of $M$ are nondegenerate, since they are
pairwise disjoint intervals; their union $A$ is an $F_\sigma$ set, so $B=M\setminus A$ is an uncountable Borel set. Let
$P\subset B$ be a Cantor-like set. If an interval $(p,q)$ contiguous to $P$ were a subset of $M$, then $[p,q]\subset M$ would
lie in a nondegenerate component, so $p\in A$, which is impossible. Let $a$ be a bilateral condensation point of $P$ and let
$K_n$, $M^n$ be as in Theorem~\ref{thm:KS6}. Each $M^n$ is closed, and $M\setminus M^n=\{a\}\cup\bigcup_{k\neq n}M^k$ is closed
because $M^k\to a$; so $M^n$ is clopen. The last property follows from Theorem~\ref{thm:KS4}.
\end{proof}

\subsection{Components}
For a compact set $M\subset\R$ we write $c(M)$ for the cardinality of its set of components and $r(M)$ for the cardinality
of its set of nondegenerate components.

\begin{lemma}\label{lem:components}
Let $\pi\colon M\to N$ be a continuous surjection between compact subsets of $\R$.
\begin{enumerate}[\upshape(a)]
\item $c(N)\le c(M)$.
\item If $c(M)\le\aleph_0$, then $r(N)\le r(M)$.
\item If $c(M)\le\aleph_0$ and $r(M)=r(N)<\infty$, then $\pi$ maps each nondegenerate component of $M$ onto a nondegenerate
component of $N$, this defines a bijection between the nondegenerate components, and $\pi$ maps the union of the degenerate
components of $M$ onto a set containing the union of the degenerate components of $N$.
\end{enumerate}
\end{lemma}

\begin{proof}
(a) The image of a component of $M$ is connected, so it lies in one component of $N$. For every component $D$ of $N$ choose
$y\in D$, a preimage $x$ of $y$ and the component $C_D$ of $x$; then $\pi(C_D)\subset D$, so $D\mapsto C_D$ is injective.

(b) Let $J$ be a nondegenerate component of $N$. It is the union of the sets $\pi(C)$ over the components $C$ of $M$ with
$\pi(C)\subset J$. If all these $C$ were degenerate, $J$ would be countable. Choose a nondegenerate $C_J$ with
$\pi(C_J)\subset J$; again $J\mapsto C_J$ is injective.

(c) By (b) and finiteness the map $J\mapsto C_J$ is a bijection, so $C_J$ is the only nondegenerate component of $M$ mapped
into $J$. Hence $J\setminus\pi(C_J)$ is contained in the image of countably many points, and since $\pi(C_J)$ is a compact
interval, $\pi(C_J)=J$. Thus $\pi$ maps the union $U_M$ of the nondegenerate components of $M$ onto the union $U_N$ of
those of $N$, and the remaining points of $N$ are images of points of $M\setminus U_M$.
\end{proof}

\section{Realization over a cycle: proof of Theorem~\ref{thm:A}}\label{sec:A}

\begin{lemma}\label{lem:phase}
Let $f\colon\I\to\I$ be continuous and let $x_0\in\I$ have a finite $\omega$-limit set. Then $\omega_f(x_0)$ is a periodic orbit
$P=\{p_0,\dots,p_{m-1}\}$, $f(p_i)=p_{i+1}$, and there is $\nu\in\Z_m$ with $f^{km+i}(x_0)\to p_{i+\nu}$ as $k\to\infty$, for every
$i$.
\end{lemma}

\begin{proof}
The set $\omega_f(x_0)$ is nonempty, finite and satisfies $f(\omega_f(x_0))=\omega_f(x_0)$, and it is \ICT\ by
Lemma~\ref{lem:pseudo-ICT} applied to the orbit. A finite \ICT\ set is a single periodic orbit, since for $\eps$ smaller
than the distances between its points every $\eps$-chain in it is a true orbit. Let $\delta>0$ be such that the sets $B(p_i,\delta)$ are pairwise disjoint, and choose $0<\eta<\delta$ with $f(B(p_i,\eta))\subset B(p_{i+1},\delta)$ for every $i$.
Since $\dist(f^n(x_0),P)\to0$, there is $N$ with $f^n(x_0)\in\bigcup_iB(p_i,\eta)$ for $n\ge N$. If $f^n(x_0)\in B(p_i,\eta)$
with $n\ge N$, then $f^{n+1}(x_0)$ lies in $B(p_{i+1},\delta)$ and in some $B(p_k,\eta)\subset B(p_k,\delta)$, so $k=i+1$.
Hence the index advances by one at each step after time $N$, which gives the claim.
\end{proof}

We next show that the hypotheses on the base in the ``moreover'' part of Theorem~\ref{thm:A} can always be met.

\begin{lemma}\label{lem:base}
For every $m\ge1$ there are a piecewise linear map $f\colon\I\to\I$, a periodic orbit $P=\{p_0,\dots,p_{m-1}\}$ with
$f(p_i)=p_{i+1}$, and a point $x_0$ with $f^{km+i}(x_0)\to p_i$ for every $i$ and $f^n(x_0)\notin P$ for all $n$.
\end{lemma}

\begin{proof}
Let $p_i=\frac{i+1}{m+1}$, $\delta=\frac1{4(m+1)}$ and $U_i=[p_i-\delta,p_i+\delta]$. The intervals $U_i$ are pairwise disjoint
and contained in $(0,1)$. On $U_i$ put $f(x)=p_{i+1}+\frac12(x-p_i)$, so $f(U_i)\subset U_{i+1}$, and extend $f$ linearly on the
gaps between the $U_i$ and constantly on $[0,p_0-\delta]$ and $[p_{m-1}+\delta,1]$. For $x_0=p_0+\delta$ we get
$f^n(x_0)=p_{n\bmod m}+2^{-n}\delta$.
\end{proof}

\begin{proof}[Proof of \textup{(a)$\Rightarrow$(b)} and of one half of the last assertion]
Let $F=(f,g)$, $z=(x_0,y_0)$ and $P$ realize $\M$. By Lemma~\ref{lem:phase}, after shifting cyclically the indices of the
points of $P$ and, simultaneously, of the sets $M_i$ (which does not affect (b)), we may assume that $x_n:=f^n(x_0)$ satisfies
$x_{km+i}\to p_i$; in the ``moreover'' part this holds by hypothesis. Write
$F^n(z)=(x_n,y_n)$, put $\varphi_i=g(p_i,\cdot)$, and let
$u_n=(n\bmod m,\,y_n)\in\Z_m\times\I$. By uniform continuity of $g$,
\[
|y_{n+1}-\varphi_{n\bmod m}(y_n)|\le\sup_{y\in\I}\bigl|g(x_n,y)-g(p_{n\bmod m},y)\bigr|\longrightarrow0 ,
\]
so $(u_n)$ is an asymptotic pseudo-orbit of the induced cyclic map $\Phi$. Since $x_{km+i}\to p_i$, the accumulation points of
$(u_n)$ in $\{i\}\times\I$ are exactly the points of $\{i\}\times M_i$. Thus $\omega((u_n))=\bigcup_i\{i\}\times M_i$, which is
\ICT\ for $\Phi$ by Lemma~\ref{lem:pseudo-ICT}. By Lemma~\ref{lem:invariant}, $\Phi$ maps this set onto itself, and its
restriction is a cyclic map of $X_\M$, which is chain transitive because chains in the set only use $\Phi$ there.
\end{proof}

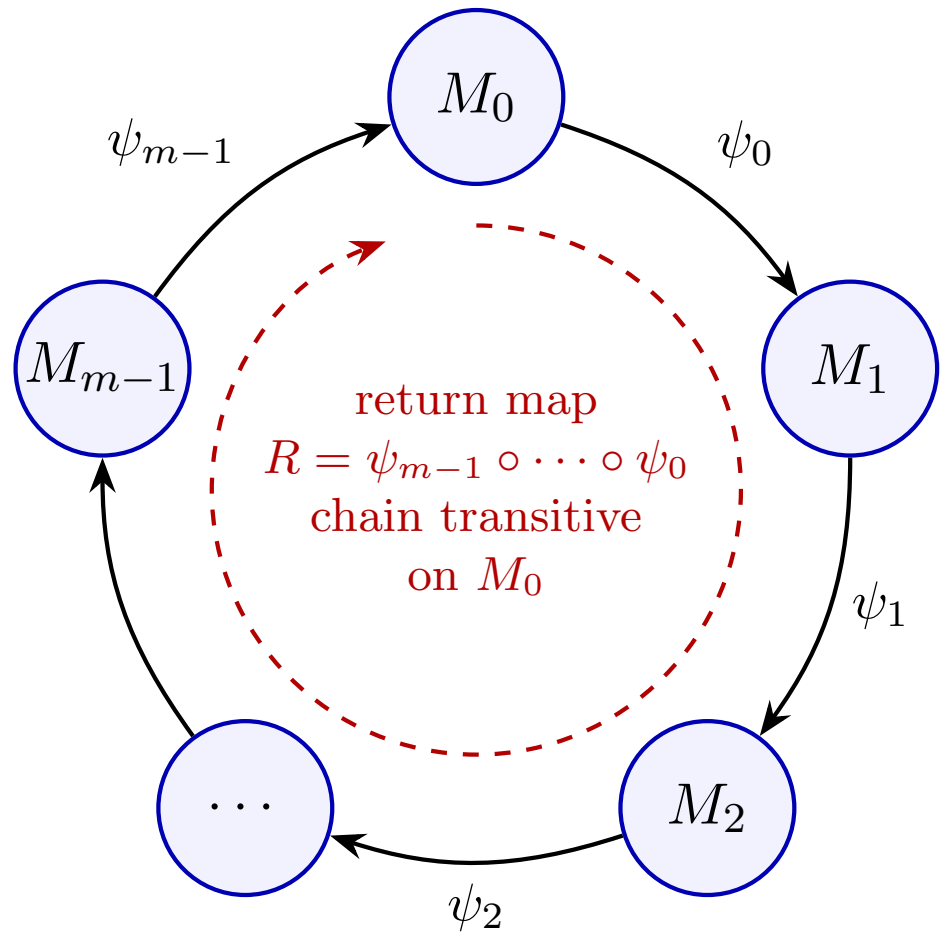
\begin{figure}[t]
\centering
\begin{tikzpicture}[scale=1,>=Stealth]
  \foreach \i/\ang in {0/90,1/18,2/-54,3/-126,4/162}{
     \node[circle,thick,minimum size=1.15cm,draw=blue!70!black,fill=blue!5,text=blue!70!black] (M\i) at (\ang:2.6) {};
  }
  \node at (M0) {$M_0$}; \node at (M1) {$M_1$}; \node at (M2) {$M_2$};
  \node at (M3) {$\cdots$}; \node at (M4) {$M_{m-1}$};
  \draw[->,thick] (M0) to[bend left=18] node[above right,font=\small]{$\psi_0$} (M1);
  \draw[->,thick] (M1) to[bend left=18] node[right,font=\small]{$\psi_1$} (M2);
  \draw[->,thick] (M2) to[bend left=18] node[below,font=\small]{$\psi_2$} (M3);
  \draw[->,thick] (M3) to[bend left=18] (M4);
  \draw[->,thick] (M4) to[bend left=18] node[above left,font=\small]{$\psi_{m-1}$} (M0);
  \draw[->,red!70!black,thick,dashed] (0,1.75) arc (90:-250:1.75);
  \node[red!70!black,font=\footnotesize,align=center] at (0,0) {return map\\ $R=\psi_{m-1}\circ\cdots\circ\psi_0$\\ chain transitive\\ on $M_0$};
\end{tikzpicture}
\caption{Theorem~\ref{thm:A}(c). A chain of the cyclic map $\Psi$ advances one section at a time; grouping its steps in
blocks of $m$ gives a chain of the return map $R$, and conversely each step of a chain of $R$ unfolds into $m$ steps of
$\Psi$, all exact but the last one.}
\label{fig:return}
\end{figure}

\begin{proof}[Proof of \textup{(b)$\Leftrightarrow$(c)}]
(b)$\Rightarrow$(c). Let $\Psi$ be a chain transitive cyclic map and $\psi_i\colon M_i\to M_{i+1}$ its components. By
Lemma~\ref{lem:invariant}, $\Psi(X_\M)=X_\M$; since only $\{i\}\times M_i$ is mapped into $\{i+1\}\times M_{i+1}$, every
$\psi_i$ is onto. Let $a,b\in M_0$ and $\eps>0$, take $\eta<1$ as in Lemma~\ref{lem:blocks} (applied to $\Psi$ with this $m$),
and an $\eta$-chain of $\Psi$ from $(0,a)$ to $(0,b)$. Since pieces are at distance $2$ and $\eta<1$, each step advances the
index by one, so the chain has length $sm$ with $s\ge1$. Taking every $m$-th point gives an $\eps$-chain of
$R=\Psi^m|_{M_0}$ from $a$ to $b$.

(c)$\Rightarrow$(b). Define $\Psi(i,y)=(i+1,\psi_i(y))$. Let $(i,a),(j,b)\in X_\M$ and $\eps>0$. Follow the exact orbit from
$(i,a)$ to $(0,a')$, $a'=\psi_{m-1}\circ\dots\circ\psi_i(a)$ (no steps if $i=0$). By surjectivity choose $c\in M_0$ with
$\psi_{j-1}\circ\dots\circ\psi_0(c)=b$ ($c=b$ if $j=0$). Take an $\eps$-chain $a'=v_0,v_1,\dots,v_s=c$ of $R$ in $M_0$ and
replace each step $v_t\to v_{t+1}$ by the $m$ steps
\[
(0,v_t)\to\Psi(0,v_t)\to\dots\to\Psi^{m-1}(0,v_t)\to(0,v_{t+1}),
\]
all exact except the last one, whose error is $|R(v_t)-v_{t+1}|<\eps$. Finally follow the exact orbit from $(0,c)$ to $(j,b)$.
The result is an $\eps$-chain of $\Psi$ from $(i,a)$ to $(j,b)$ of positive length (Figure~\ref{fig:return}).
\end{proof}

\begin{proof}[Proof of \textup{(b)$\Rightarrow$(a)} and of the other half of the last assertion]
Let $f$, $P$, $x_0$ and $\varphi_i$ be as in the ``moreover'' statement (they exist by Lemma~\ref{lem:base}), and let
$K=\bigcup_i\{i\}\times M_i$ be closed and \ICT\ for the induced map $\Phi$.
If we only assume (b), let $\Psi$ be a chain transitive cyclic map with components $\psi_i$. Extend each $\psi_i$ to a continuous
$\varphi_i\colon\I\to\I$, interpolating linearly on the bounded components of $\I\setminus M_i$ and constantly on the others;
then $K$ is \ICT\ for the induced map, because chains in $K$ only see $\Phi|_K=\Psi$. So it suffices to treat the
``moreover'' situation. The argument follows \cite[Lemma~2]{KS}.

\emph{Step 1: the fibre sequence.} Apply Lemma~\ref{lem:ICT-pseudo} to $\Phi$ and $K$, starting at a point of
$\{0\}\times M_0$. We get an asymptotic pseudo-orbit $u_n=(i_n,y_n)$ in $K$ with $\omega((u_n))=K$ and all one-step errors
$<1$. Since different pieces are at distance $2$, $i_{n+1}=i_n+1$, so $i_n=n\bmod m$. Let
$\delta_n=|y_{n+1}-\varphi_{i_n}(y_n)|\to0$ and
\[
h_n(y)=\min\bigl\{1,\max\{0,\ \varphi_{i_n}(y)+y_{n+1}-\varphi_{i_n}(y_n)\}\bigr\},\qquad y\in\I .
\]
Then $h_n$ is continuous, $h_n(y_n)=y_{n+1}$ and $\|h_n-\varphi_{i_n}\|_\infty\le\delta_n$.

\emph{Step 2: the fibre maps.} The points $x_n=f^n(x_0)$ are pairwise distinct: a repetition would make the orbit
eventually periodic, and since it accumulates on $P$ it would have to enter $P$, which is excluded. The set
$E=P\cup\{x_n:n\ge0\}$ is closed, and every $x_n$ is isolated in $E$. Define $G$ on $E\times\I$ by $G(x_n,y)=h_n(y)$ and
$G(p_i,y)=\varphi_i(y)$. Then $G$ is continuous: at points $(x_n,y)$ because $x_n$ is isolated, and at $(p_i,y)$ because
$(x_{n_k},y_k)\to(p_i,y)$ forces $n_k\to\infty$ and $n_k\equiv i\pmod m$ for large $k$, so that
$|h_{n_k}(y_k)-\varphi_i(y)|\le\delta_{n_k}+|\varphi_i(y_k)-\varphi_i(y)|\to0$. By the Tietze extension theorem $G$ extends to a
continuous $g\colon\I^2\to\I$ with $g(p_i,\cdot)=\varphi_i$.

\emph{Step 3: the orbit.} For $F=(f,g)$ and $z=(x_0,y_0)$, induction gives $F^n(z)=(x_n,y_n)$. Since $x_{km+i}\to p_i$ and the
accumulation points of $(y_{km+i})_k$ form $M_i$, we obtain $\omega_F(z)=\bigcup_i\{p_i\}\times M_i$.
\end{proof}

\begin{proof}[Proof of Corollary~\ref{cor:B}]
(ii) is part of Theorem~\ref{thm:A}(c). For (i), let $\Psi$ be a chain transitive cyclic map of $X_\M$. The proof of
(b)$\Rightarrow$(c) applies to every piece, so $\Psi^m$ restricted to $\{i\}\times M_i$ is chain transitive. Hence the
$1$-tuple $(M_i)$ satisfies (c) of Theorem~\ref{thm:A} with $m=1$, so it is realizable, and $M_i\in\Fib$ by
Theorem~\ref{thm:KS}.
\end{proof}

\begin{remark}\label{rem:KS}
For $m=1$, Theorem~\ref{thm:A} and Theorem~\ref{thm:KS} show that a nonempty compact set $K\subset\I$ admits a chain transitive
continuous self-map if and only if $K\in\Fib$; by an affine change of coordinates (Lemma~\ref{lem:metric}) the same holds for
compact subsets of $\R$. This is essentially contained in \cite{KS}. There a map $h\colon K\to K$ is said to belong to
$K(\mathcal E)$ if for every $\eps>0$ some cyclic $\eps$-chain of $h$ is an $\eps$-net of $K$. A routine argument with
uniform continuity shows that this is equivalent to chain transitivity of $h$. The proof of \cite[Theorem~1]{KS} constructs such a map for
every $K\in\Fib$ that is neither nowhere dense nor a finite union of intervals; for these two classes it refers to
Theorem~\ref{thm:ABCP}, from which such a map is obtained by Lemma~\ref{lem:pseudo-ICT}. The converse part of the proof uses
only that the map is onto and weakly incompressible, which chain transitive maps are \cite{BGOR}.
\end{remark}

\begin{remark}\label{rem:homeo}
If $M_0\in\Fib$ and $h_i\colon M_0\to M_i$ are homeomorphisms, with $h_0=\mathrm{id}$, then $\M$ is realizable. Indeed, take a
chain transitive $\rho$ on $M_0$ (Remark~\ref{rem:KS}), $\psi_i=h_{i+1}\circ h_i^{-1}$ for $i\le m-2$ and
$\psi_{m-1}=\rho\circ h_{m-1}^{-1}$. Then $R=\rho$, and Theorem~\ref{thm:A}(c) applies. This recovers the Kolyada--Snoha result
on products $B\times M$ when $B$ is a periodic orbit.
\end{remark}

\section{Uncountably many components: proof of Theorem~\ref{thm:C}}\label{sec:C}

\begin{proof}[Proof of Theorem~\ref{thm:C}]
For each $i\in\Z_m$ apply Lemma~\ref{lem:uncountable} to $M_i$. We get $a_i\in M_i$ and pairwise disjoint nonempty clopen sets
$M_i^n$, $n\ge0$, with $M_i=\{a_i\}\cup\bigcup_nM_i^n$, $M_i^n\to a_i$, and continuous surjections from each $M_i^n$ onto any
nonempty compact subset of $\I$. Define $\psi_i\colon M_i\to M_{i+1}$ by $\psi_i(a_i)=a_{i+1}$ and
\begin{align*}
&\psi_i \text{ maps } M_i^n \text{ onto } M_{i+1}^n &&\text{for } 0\le i\le m-2,\ n\ge0,\\
&\psi_{m-1} \text{ maps } M_{m-1}^n \text{ onto } M_0^{n-1} &&\text{for } n\ge1,\\
&\psi_{m-1}(M_{m-1}^0)=\{a_0\}.
\end{align*}
Each $\psi_i$ is continuous on the clopen sets $M_i^n$. It is continuous at $a_i$: given a neighbourhood $V$ of $a_{i+1}$,
there is $N$ with $M_{i+1}^{n}\subset V$ for $n\ge N-1$, and a neighbourhood of $a_i$ that meets only $a_i$ and sets
$M_i^n$ with $n\ge N$, because the finitely many sets $M_i^n$ with $n<N$ are closed and do not contain $a_i$. Each $\psi_i$
is onto. Let $\Psi$ be the cyclic map with components $\psi_i$.

Every point of $X_\M$ reaches $(0,a_0)$ by an exact orbit: a point of $\{i\}\times M_i^n$ reaches
$\{m-1\}\times M_{m-1}^n$, then $\{0\}\times M_0^{n-1}$, and so on, and after finitely many steps lands on $(0,a_0)$; the points
$(i,a_i)$ reach it as well. Conversely, let $(j,b)\in X_\M$ and $\eps>0$. Choose $N$ so large that $M_0^n$ lies in the
$\eps$-neighbourhood of $a_0$ for $n\ge N$ and that $b\in\{a_j\}\cup\bigcup_{n<N}M_j^n$. The orbit of $(0,a_0)$ reaches
$(m-1,a_{m-1})$, whose image is $(0,a_0)$. If $b\neq a_j$, say $b\in M_j^n$ with $n<N$, then the composition of the maps
$M_0^N\to M_1^N\to\dots\to M_{m-1}^N\to M_0^{N-1}\to\dots\to M_j^n$, each of which is onto, is onto $M_j^n$; choose $c\in M_0^N$
whose image under it is $b$, and replace the image $(0,a_0)$ of $(m-1,a_{m-1})$ by $(0,c)$, an error smaller than $\eps$. Then an
iterate of $(0,c)$ equals $(j,b)$. If $b=a_j$, follow the exact orbit of $(0,a_0)$. Concatenating, every point is joined to every point by an $\eps$-chain of positive length, so $\Psi$ is chain
transitive, and $\M$ is realizable by Theorem~\ref{thm:A}.
\end{proof}

\section{Countably many components}\label{sec:D}

\subsection{Countable sections}

\begin{proposition}\label{prop:countable}
Let $\M$ be an $m$-tuple of nonempty countable closed subsets of $\I$. The following are equivalent: $\M$ is realizable;
$\M$ satisfies \textup{(ii)} of Corollary~\ref{cor:B}; the sets $M_i$ are pairwise homeomorphic.
\end{proposition}

\begin{proof}
Realizability implies (ii) by Corollary~\ref{cor:B}. Assume (ii). By Lemma~\ref{lem:CB} the ranks of the $M_i$ do not
increase along the cycle of surjections, so they are all equal, and then the numbers of points of the kernels do not
increase either, so they are all equal. By Theorem~\ref{thm:MS} the $M_i$ are pairwise homeomorphic. Finally, a countable
set is never KS-exceptional, since Definition~\ref{def:exceptional} requires a nondegenerate interval; so homeomorphic
countable sections give a realizable tuple by Remark~\ref{rem:homeo}.
\end{proof}

\subsection{Finitely many nondegenerate components}
In this subsection $\M$ is an $m$-tuple satisfying (i) and (ii) such that every $M_i$ has countably many components and
exactly $r$ nondegenerate components, where $1\le r<\infty$; by Lemma~\ref{lem:components}(b) the number $r$ is the same
for all $i$. Write $M_i=U_i\cup C_i$, where $U_i$ is the union of the nondegenerate components and $C_i$ is the union of the
degenerate ones. We fix continuous surjections $\sigma_i\colon M_i\to M_{i+1}$, $i\in\Z_m$.

By Lemma~\ref{lem:components}(c), $\sigma_i$ maps each nondegenerate component $J$ of $M_i$ onto a nondegenerate component
$\tau_i(J)$ of $M_{i+1}$, the map $\tau_i$ is a bijection, and $C_{i+1}\subset\sigma_i(C_i)$. In particular, if some $C_i$ is
empty, then all of them are empty, every $M_i$ is a union of $r$ disjoint nondegenerate compact intervals, the $M_i$ are
pairwise homeomorphic, and $\M$ is realizable by Remark~\ref{rem:homeo}, because such a set is not KS-exceptional.

So assume that every $C_i$ is nonempty. By (i), $\dist(C_i,J)=0$ for every nondegenerate component $J$ of $M_i$; since a point
of $\I$ close to an interior point of $J$ and belonging to $M_i$ lies in $J$, $C_i$ accumulates at an endpoint of $J$. In
particular $C_i$ is infinite. Let $K_i=\overline{C_i}$, a countable compact set, and let $E_i=K_i\cap U_i$, the set of
endpoints of nondegenerate components of $M_i$ at which $C_i$ accumulates. Since $K_i$ is infinite, $\rk K_i\ge1$.

\begin{lemma}\label{lem:tops}
There are an ordinal $\alpha\ge1$ and an integer $n\ge1$ such that $\rk K_i=\alpha$ and $|\Ker K_i|=n$ for every $i$. Put
$T_i=\Ker K_i$. Then $\sigma_i$ maps $T_i$ bijectively onto $T_{i+1}$ and $T_i\cap E_i$ onto $T_{i+1}\cap E_{i+1}$, and if
$t\in T_i$ lies in a nondegenerate component $J$, then $\sigma_i(t)$ is an endpoint of $\tau_i(J)$. Consequently $J$ and
$\tau_i(J)$ contain the same number of points of $T_i$ and $T_{i+1}$.
\end{lemma}

\begin{proof}
Since $C_{i+1}\subset\sigma_i(C_i)$ and $\sigma_i(K_i)$ is compact, $K_{i+1}\subset\sigma_i(K_i)$. Let
$Z_i=\sigma_i^{-1}(K_{i+1})\cap K_i$; it is compact and $\sigma_i(Z_i)=K_{i+1}$. By Lemma~\ref{lem:CB},
$K_{i+1}^{(\beta)}\subset\sigma_i(Z_i^{(\beta)})\subset\sigma_i(K_i^{(\beta)})$ for every $\beta$. Hence $\rk K_{i+1}\le\rk K_i$,
so all ranks are equal to some $\alpha$, and then $|T_{i+1}|\le|T_i|$, so all kernels have the same number $n$ of points and
$\sigma_i(T_i)\supset T_{i+1}$ forces $\sigma_i(T_i)=T_{i+1}$ with $\sigma_i$ injective on $T_i$. If $t\in T_i\cap J$ with $J$
a nondegenerate component, then $\sigma_i(t)\in\sigma_i(J)=\tau_i(J)$ and $\sigma_i(t)\in K_{i+1}$, and
$K_{i+1}\cap\tau_i(J)$ consists of endpoints of $\tau_i(J)$. So $\sigma_i$ maps $T_i\cap E_i$ injectively into
$T_{i+1}\cap E_{i+1}$; along the cycle these finite numbers cannot increase, so they are equal and the map is onto. The last
assertion follows because $\sigma_i$ is injective on $T_i$ and maps $T_i\cap J$ into $\tau_i(J)$, while every point of
$T_{i+1}\cap\tau_i(J)$ is the image of a point of $T_i$ lying in a component mapped into $\tau_i(J)$, that is, in $J$.
\end{proof}

\begin{proposition}\label{prop:finite-r}
Under the assumptions of this subsection, with every $C_i$ nonempty, $\M$ is realizable.
\end{proposition}

\begin{proof}
We construct a chain transitive cyclic map $\Psi$ of $X_\M$ and apply Theorem~\ref{thm:A}. We drop the index $i$ from
$\Psi$ and the first coordinate from points when no confusion arises, and we write $\sigma$ for the map $\sigma_i$ on $M_i$.
Indices of sections are in $\Z_m$. Let $\mathcal J$ be the set of all nondegenerate components of all sections; the maps
$\tau_i$ define a permutation $\tau$ of the finite set $\mathcal J$. Let $T=\bigcup_iT_i$ and $T^{\mathrm f}_i=T_i\cap C_i$ (the
\emph{free tops}). When no section is specified we write $C=\bigcup_iC_i$, $E=\bigcup_iE_i$, $U=\bigcup_iU_i$, and similarly
for the sets $\mathcal N_i$, $L_i$, $T^{\mathrm f}_i$ defined below, all viewed as subsets of $X_\M$.

\emph{Step 1: maps on the intervals.} For every $J\in\mathcal J$ we choose an affine homeomorphism $h_J\colon J\to\tau(J)$
and, if $J\cap T=\emptyset$ or $|J\cap T|=1$, an endpoint $\hat e(J)\in E$ of $J$ (the \emph{designated endpoint}), such that:
\begin{enumerate}[\upshape(H1)]
\item $h_J(t)=\sigma(t)$ for every $t\in J\cap T$;
\item $h_J(\hat e(J))=\hat e(\tau(J))$ whenever $\hat e(J)$ is defined, and $\hat e(J)\in T$ if $J\cap T\neq\emptyset$;
\item for each $\tau$-cycle $(J_0,\dots,J_{\ell-1})$, the return map
$H=h_{J_{\ell-1}}\circ\dots\circ h_{J_0}\colon J_0\to J_0$ is the identity or the affine reflection of $J_0$.
\end{enumerate}
By Lemma~\ref{lem:tops} the number $|J\cap T|\in\{0,1,2\}$ is constant along a $\tau$-cycle. If it is $0$, choose
$\hat e(J_0)$ to be an endpoint of $J_0$ at which $C$ accumulates, let $\hat e(J_k)$ be such an endpoint of $J_k$ for every $k$, and let
$h_{J_k}$ be the affine map sending $\hat e(J_k)$ to $\hat e(J_{k+1})$; then $H$ is an affine homeomorphism fixing an endpoint, hence
the identity. If it is $1$, let $\hat e(J)$ be the unique point of $J\cap T$; by Lemma~\ref{lem:tops}, $\sigma(\hat e(J))=\hat e(\tau(J))$,
and the affine maps with $h_J(\hat e(J))=\hat e(\tau(J))$ satisfy (H1)--(H3), with $H$ the identity. If it is $2$, both endpoints of
$J$ lie in $T$ and $\sigma$ maps them onto the two endpoints of $\tau(J)$; let $h_J$ be the affine map that agrees with
$\sigma$ on the endpoints. Then $H$ permutes the endpoints of $J_0$, so it is the identity or the reflection.

\emph{Step 2: nodes and blocks.} Let $\mathcal N_i=T_i\cup E_i$, a finite set of \emph{nodes}, and $L_i=\mathcal N_i\setminus T_i$. For
$v\in\mathcal N_i$ let $\rho(v)=\rk_{K_i}(v)$; then $\rho(v)\ge1$ because $v$ is an accumulation point of $C_i$, and $\rho(v)<\alpha$
for $v\in L_i$. Choose pairwise disjoint clopen subsets $B_v$ of $K_i$, $v\in\mathcal N_i$, covering $K_i$, with
$B_v\cap\mathcal N_i=\{v\}$ and with $B_v^{(\rho(v))}=\{v\}$ for $v\in L_i$. This is possible: take pairwise disjoint clopen
neighbourhoods $W_v$ of the nodes with $W_v\cap K_i^{(\rho(v))}=\{v\}$ and $W_v\cap\mathcal N_i=\{v\}$, and add the clopen set $K_i\setminus\bigcup_vW_v$ to
$W_t$ for one $t\in T_i$; it contains no point of rank $\alpha$. Then $\rk B_v=\rho(v)$ and $\Ker B_v=\{v\}$ for every $v$.
Every compact subset $Q$ of $C_i$ that is clopen in $K_i$ is clopen in $M_i$: it is at positive distance from the closed set
$U_i$, so an open set $O\subset\R$ with $O\cap K_i=Q$ and $O\cap U_i=\emptyset$ satisfies $O\cap M_i=Q$.

\emph{Step 3: portions.} For every node $v$ we choose sets $V_0(v)=B_v\supset V_1(v)\supset\cdots$, clopen in $K_i$, with
$\bigcap_kV_k(v)=\{v\}$, and put $\Pi_k(v)=V_k(v)\setminus V_{k+1}(v)$, a compact subset of $C_i$ clopen in $M_i$. For
$v\in L_i$ we only require $\Pi_k(v)\neq\emptyset$ for all $k$, which is possible because $v$ is an accumulation point of
$K_i$. For the tops we require:
\begin{enumerate}[\upshape(P1)]
\item for every $k\ge1$ and all $t,t'\in T$, the set $\Pi_k(t)$ dominates $\Pi_{k-1}(t')$;
\item there are $t_0^*\in T_0$ and, writing $t_i^*$ for the image of $t_0^*$ under the composition
$\sigma_{i-1}\circ\dots\circ\sigma_0$ ($0\le i<m$), the set $\Pi_0(t_i^*)$ contains pairwise disjoint clopen subsets
$S_w$ ($w\in L_{i+1}$) with $\rk S_w\ge\rho(w)$, and nonempty pairwise disjoint clopen subsets $D_J$
($J\in\mathcal J$ a component of $M_{i+1}$), $D_t$ ($t\in T^{\mathrm f}_{i+1}$) and $O_i$, all disjoint from the $S_w$.
\end{enumerate}
We call $t_i^*$ the \emph{master} of $M_i$. Here $t^*_m$ means $t^*_0$, and $\sigma_{m-1}(t^*_{m-1})$ need not equal $t^*_0$.

The sets are chosen by recursion on $k$. For $k=0$ we choose $V_1(t)$ for all $t\in T$ so that $\Pi_0(t)\neq\emptyset$ and
(P2) holds. This is possible because $\rho(w)<\alpha$ for $w\in L_{i+1}$, so the punctured neighbourhoods of $t_i^*$ contain
infinitely many points of rank $\ge\max_w\rho(w)$ if $L_{i+1}\neq\emptyset$ (if $\alpha=\beta+1$ these are points of rank
$\beta$, which accumulate at $t_i^*$; if $\alpha$ is a limit ordinal, there are points of rank larger than
$\max_w\rho(w)$ arbitrarily close to $t_i^*$), and in any case infinitely many points; around finitely many of them we
take small pairwise disjoint clopen sets. Given $V_k(t)$ for all $t\in T$ ($k\ge1$), we choose $V_{k+1}(t)\subset V_k(t)$
of diameter $<1/k$ such that $\Pi_k(t)$ dominates the finitely many sets $\Pi_{k-1}(t')$, $t'\in T$, all of which have rank
$<\alpha$. If $\alpha=\beta+1$, all these sets have rank $\le\beta$, and the points of rank $\beta$ accumulate at $t$; we
take $V_{k+1}(t)$ so small that $\Pi_k(t)$ contains at least $\max_{t'}|\Ker \Pi_{k-1}(t')|$ such points. If $\alpha$ is a limit
ordinal, let $\gamma=\max_{t'}\rk \Pi_{k-1}(t')<\alpha$; since $t\in K_i^{(\gamma+2)}$ there is a point $y\neq t$ of rank
$\ge\gamma+1$ in $V_k(t)$, and we take $V_{k+1}(t)$ not containing $y$.

\emph{Step 4: definition of $\Psi$ on $M_i$.} Let $\mu_J$ denote the midpoint of $J\in\mathcal J$.
\begin{enumerate}[\upshape(1)]
\item On each $J\in\mathcal J$ contained in $M_i$, $\Psi=h_J$.
\item $\Psi(t)=\sigma(t)$ for $t\in T^{\mathrm f}_i$.
\item For $t\in T_i$ and $k\ge0$, $\Psi$ maps $\Pi_{k+1}(t)$ onto $\Pi_k(\sigma(t))$, by Lemma~\ref{lem:dominate} and (P1).
\item For $t\in T_i\setminus\{t_i^*\}$, $\Psi(\Pi_0(t))=\{t^*_{i+1}\}$.
\item On $\Pi_0(t^*_i)$: $\Psi$ maps $S_w$ onto $B_w$ for $w\in L_{i+1}$, by Theorem~\ref{thm:KS8}; $\Psi(D_J)=\{\mu_J\}$;
$\Psi(D_t)=\{t\}$; and $\Psi(\Pi_0(t^*_i)\setminus(\bigcup_wS_w\cup\bigcup_JD_J\cup\bigcup_tD_t))=\{t^*_{i+1}\}$.
\item Let $v\in L_i$, let $J$ be the component containing $v$, and $w=h_J(v)$.
 \begin{enumerate}[\upshape(a)]
 \item If $v=\hat e(J)$, then $J\cap T=\emptyset$, so $\tau(J)\cap T=\emptyset$ by Lemma~\ref{lem:tops} and $w=\hat e(\tau(J))\in L_{i+1}$.
 Let $\Psi$ map each $\Pi_{k+1}(v)$ to one point of $\Pi_k(w)$, and $\Psi(\Pi_0(v))=\{t^*_{i+1}\}$.
 \item If $v\neq \hat e(J)$, choose points $w_k\in\tau(J)$ with $|w_k-w|<1/(k+1)$ and let $\Psi(\Pi_k(v))=\{w_k\}$.
 \end{enumerate}
\end{enumerate}
The sets $J\in\mathcal J$, the free tops and the sets $\Pi_k(v)$ partition $M_i$, so $\Psi$ is well defined; by (H1) it agrees
with $\sigma$ on $T_i$.

\emph{Step 5: continuity and surjectivity.} $\Psi$ is continuous on each $J$ and on each clopen set $\Pi_k(v)$. Let $v$ be a
node or a free top. Since $B_v$ is a neighbourhood of $v$ in $K_i$ and $U_i$ meets a small neighbourhood of $v$ only in the
component of $v$ (if any), the points of $M_i$ near $v$ belong to that component or to sets $\Pi_k(v)$ with $k$ large. On the
component $\Psi$ is continuous, and $\Psi(\Pi_k(v))$ converges to $\Psi(v)$: in case (3) because
$\Psi(\Pi_{k+1}(t))=\Pi_k(\sigma(t))\subset V_k(\sigma(t))$ and these sets shrink to $\sigma(t)=\Psi(t)$; in case (6a) because
$\Psi(\Pi_{k+1}(v))\subset V_k(w)$ and $w=h_J(v)=\Psi(v)$; in case (6b) because $w_k\to w=\Psi(v)$. Thus $\Psi$ is continuous.
It is onto $M_{i+1}$: the maps $h_J$ cover $U_{i+1}$; every $w\in T_{i+1}$ equals $\sigma(t)=\Psi(t)$ with $t\in T_i$, and
$\Pi_k(w)=\Psi(\Pi_{k+1}(t))$; every $B_w$ with $w\in L_{i+1}$ is $\Psi(S_w)$; and the sets $B_w$, $w\in\mathcal N_{i+1}$, cover
$K_{i+1}\supset C_{i+1}$.

\emph{Step 6: chain transitivity.} Fix $\eps>0$ and let $Q=\bigcup_i(U_i\cup T^{\mathrm f}_i)$.

\emph{Claim 1. For every $x\in X_\M$ some $\Psi^k(x)$, $k\ge0$, lies in $Q$.} If $x\in \Pi_k(t)$ with $t\in T$, then by (3)
$\Psi^k(x)$ lies in some $\Pi_0(t')$, and by (4) and (5) $\Psi^{k+1}(x)$ is a master, a midpoint, a free top, or a point of some
$B_w$ with $w\in L$. A master lies in $Q$. If $x\in \Pi_k(v)$ with $v\in L$, then in case (6b) $\Psi(x)\in U$, and in case (6a)
the orbit stays in sets $\Pi_j(w')$ with $w'$ designated endpoints in $L$ and $j$ decreasing, and after $k+1$ steps it is a
master. Nodes in $L$ lie in $U$.

\emph{Claim 2. Every $\tau$-cycle is internally chain connected.} If $x,y$ lie in the union of the components of a
$\tau$-cycle $(J_0,\dots,J_{\ell-1})$, there is an $\eps$-chain from $x$ to $y$ inside this union. Indeed, follow the exact
orbit of $x$ to a point $x'\in J_0$, take an $\eps$-chain of $H$ in $J_0$ from $x'$ to the point $y'$ whose orbit reaches $y$
(Lemma~\ref{lem:interval-CT}), unfold each step of it into $\ell$ steps of $\Psi$ as in the proof of
(c)$\Rightarrow$(b) of Theorem~\ref{thm:A}, and follow the exact orbit of $y'$ to $y$.

\emph{Claim 3. From every $x\in Q$ there is an $\eps$-chain to the master $t^*_0$.} Let first $x\in J$ with
$J\in\mathcal J$. The $\tau$-cycle of $J$ contains a point $q$ that is either a designated endpoint or an endpoint in $T$; by
Claim 2 there is an $\eps$-chain from $x$ to $q$. The point $q'=\Psi(q)$ is a node and $\Pi_k(q')\to q'$, so there is $N\ge1$
such that $\Pi_N(q')$ lies in the $\eps$-neighbourhood of $q'$; we shall replace the image of $q$ by a suitable point
$c\in\Pi_N(q')$. If $q'\in L$, every such $c$ reaches a master by (6a). If $q'\in T$, the composition of the maps in (3) sends
$\Pi_N(q')$ onto $\Pi_0(t)$ for some top $t$, so we may choose $c$ whose $N$-th image lies in the set mapped to a master by
(4) or (5); when $t$ is the master this set contains the nonempty set $O_i$. If $x\in T^{\mathrm f}_i$, replace $\Psi(x)\in T$
by a suitable point $c\in\Pi_N(\Psi(x))$ in the same way. In all
cases we reach some master $t^*_j$, and since $\Psi$ permutes the finite set $T$, the exact orbit of $t^*_j$ reaches
$t^*_0$, which lies in the same $\Psi$-orbit.

\emph{Claim 4. Every $y\in X_\M$ can be reached from $t^*_0$ by an $\eps$-chain.} Say that $y$ is \emph{reachable} if
there is an $\eps$-chain from $t^*_0$ to $y$. Since $\Psi$ permutes $T$, the exact orbit of $t^*_0$ returns to $t^*_0$, so
$t^*_0$ and every point of its $\Psi$-orbit are reachable. Images of reachable points are reachable. Let $s\in T$ be reachable,
let $p\ge1$ be its period under $\Psi$, and let $y\in \Pi_k(\Psi^j(s))$ with $j\ge1$. Choose $\theta$ so large that, with
$N=k+j-1+\theta p$, the set $\Pi_N(\Psi(s))$ lies in the $\eps$-neighbourhood of $\Psi(s)$, and replace the image of $s$ by a
suitable point $c\in\Pi_N(\Psi(s))$. By (3), $\Psi^{j-1+\theta p}$ maps
$\Pi_N(\Psi(s))$ onto $\Pi_k(\Psi^{j+\theta p}(s))=\Pi_k(\Psi^j(s))$, so $y$ is reachable. Applying this to $s=t^*_0$ shows that all
points of $\Pi_0(t^*_i)$ are reachable for every $i$ (the masters lie on the orbit of $t^*_0$), hence, by (5), every point of
every $B_w$ with $w\in L$, every midpoint $\mu_J$ and every free top. By Claim 2 every point of every $J\in\mathcal J$ is
reachable, in particular every top lying in $U$. So every top is reachable, and applying the first part to all tops, every
point of every $\Pi_k(t)$, $t\in T$, is reachable. These sets, the sets $B_w$, the components in $\mathcal J$ and the free
tops cover $X_\M$.

By Claims 1, 3 and 4, any two points of $X_\M$ are joined by an $\eps$-chain of positive length. So $\Psi$ is chain
transitive and $\M$ is realizable by Theorem~\ref{thm:A}.
\end{proof}

\subsection{Proofs of Theorems~\ref{thm:D} and~\ref{thm:F} and of Corollary~\ref{cor:E}}

\begin{proof}[Proof of Theorem~\ref{thm:D}]
Necessity is Corollary~\ref{cor:B}. Assume (i) and (ii). By Lemma~\ref{lem:components}(b), applied along the cycle of
surjections given by (ii), all $M_i$ have the same number $r<\infty$ of nondegenerate components. If $r=0$, the sets are countable and Proposition~\ref{prop:countable} applies. If $r\ge1$, the
sets satisfy the assumptions of the previous subsection, where realizability was shown.
\end{proof}

\begin{proof}[Proof of Theorem~\ref{thm:F}]
By Lemma~\ref{lem:components}(a) the cardinalities $c(M_i)$ do not increase along the cycle of surjections given by (ii), so
either all of them are uncountable or all of them are countable. In the second case Lemma~\ref{lem:components}(b) shows that
the numbers $r(M_i)$ do not increase along the cycle, so they are all equal. If the $M_i$ have uncountably many components,
$\M$ is realizable by Theorem~\ref{thm:C}. If they have countably many components and $r(M_i)<\infty$, Conjecture~\ref{conj}
holds by Theorem~\ref{thm:D}. The remaining case is $M_i\in\W$ for all $i$. A set in $\W$ has infinitely many nondegenerate
components, so it is not KS-exceptional and (i) holds automatically.
\end{proof}

\begin{proof}[Proof of Corollary~\ref{cor:E}]
If the $M_i$ are pairwise homeomorphic and each is finite or a finite union of nondegenerate intervals, then $M_0\in\Fib$ and
$\M$ is realizable by Remark~\ref{rem:homeo}. Conversely, let $\M$ be realizable. By Corollary~\ref{cor:B} and
Lemma~\ref{lem:components}, all $M_i$ have the same number $k$ of components and the same number $r$ of nondegenerate ones.
If $0<r<k$, each $M_i$ is a finite union of nondegenerate intervals together with a nonempty finite set at positive distance
from them, which is KS-exceptional, contradicting (i). Hence either $r=0$ and every $M_i$ consists of $k$ points, or $r=k$
and every $M_i$ is a union of $k$ disjoint nondegenerate intervals; in both cases the $M_i$ are pairwise homeomorphic.
\end{proof}

\section{Examples}\label{sec:examples}

A convenient base for $m=2$ is $f(x)=(1-x)^2$. It has the periodic orbit $\{0,1\}$ and maps $(0,1)$ into $(0,1)$. Since
$f^2(x)=x^2(2-x)^2\le4x^2\le x/2$ for $0\le x\le1/8$, the orbit of $x_0=1/8$ satisfies $f^{2k}(x_0)\to0$ and
$f^{2k+1}(x_0)\to1$ without reaching the orbit $\{0,1\}$. So $p_0=0$, $p_1=1$, and the hypotheses on the base in the
``moreover'' part of Theorem~\ref{thm:A} hold. Let $\mathcal C$ denote the middle-thirds Cantor set.

\begin{example}\label{ex:cantor}
Let $M_0=\mathcal C$ and $M_1=\frac13\mathcal C\cup[\frac23,1]$. Both sets have uncountably many components, so $(M_0,M_1)$
is realizable by Theorem~\ref{thm:C}: there is a triangular map $F(x,y)=((1-x)^2,g(x,y))$ with
\[
\omega_F\bigl(\tfrac18,y_0\bigr)=\{0\}\times\mathcal C\ \cup\ \{1\}\times\bigl(\tfrac13\mathcal C\cup[\tfrac23,1]\bigr).
\]
One section is nowhere dense and the other one has nonempty interior.
\end{example}

\begin{example}\label{ex:nonhomeo}
Let $J=[\frac12,\frac34]$. For $k\ge1$ and $l\ge0$ let $b_k=\frac12-2^{-k-1}$ and $a_{k,l}=b_k-2^{-k-l-3}$, and let
$\Lambda=\{b_k\}_{k\ge1}\cup\{a_{k,l}\}_{k\ge1,\,l\ge0}$. The set $\Lambda\cup\{\frac12\}$ is closed, countable and of rank $2$ with
kernel $\{\frac12\}$; the points $b_k$ have rank $1$. Let $S=\{\frac34+2^{-k-3}:k\ge0\}$ and
\[
M_0=\Lambda\cup J,\qquad M_1=\Lambda\cup J\cup S .
\]
Both sets have countably many components, exactly one of which is nondegenerate, and in both the degenerate components
accumulate at $J$, so (i) holds. For (ii), let $W_k=\{b_k\}\cup\{a_{k,l}:l\ge0\}$; each $W_k$ is clopen in $M_0$ and
homeomorphic to a convergent sequence with its limit. A continuous surjection $M_0\to M_1$ is given by the identity on $J$,
a homeomorphism of $W_1$ onto $S\cup\{\frac34\}$ sending $b_1$ to $\frac34$, and homeomorphisms $W_{k+1}\to W_k$ for
$k\ge1$. A continuous surjection $M_1\to M_0$ is given by the identity on $M_0$ and $\frac34+2^{-k-3}\mapsto\frac34-2^{-k-3}$
on $S$. By Theorem~\ref{thm:D} the tuple $(M_0,M_1)$ is realizable. The sets are not homeomorphic: in $M_1$ both endpoints of
the unique nondegenerate component are limits of other components, in $M_0$ only one is.
\end{example}

\begin{example}\label{ex:necessary}
The tuples $([0,1],\,[0,\frac13]\cup[\frac23,1])$ and $(\mathcal C,[0,1])$ are not realizable over a periodic orbit of
period $2$, because a connected set has no continuous surjection onto a disconnected one, so (ii) fails. No tuple containing
$[0,\frac12]\cup\{1\}$ is realizable, because this set is KS-exceptional, so (i) fails.
\end{example}

\begin{table}[t]
\centering
\renewcommand{\arraystretch}{1.25}
\begin{tabular}{@{}lll@{}}
\toprule
Tuple $(M_0,M_1)$ & Realizable? & Reason\\
\midrule
$(\mathcal C,\ \frac13\mathcal C\cup[\frac23,1])$ & yes & Theorem~\ref{thm:C}\\
$(\Lambda\cup J,\ \Lambda\cup J\cup S)$ of Example~\ref{ex:nonhomeo} & yes & Theorem~\ref{thm:D}\\
$(K,K')$ with $K\cong K'\in\Fib$ & yes & Remark~\ref{rem:homeo}\\
$([0,1],\ [0,\frac13]\cup[\frac23,1])$ & no & Corollary~\ref{cor:B}(ii)\\
$(\mathcal C,\ [0,1])$ & no & Corollary~\ref{cor:B}(ii)\\
$([0,\frac12]\cup\{1\},\ M_1)$ & no & Corollary~\ref{cor:B}(i)\\
\bottomrule
\end{tabular}
\caption{Realizability of some tuples over a periodic orbit of period $2$.}
\label{tab:examples}
\end{table}

\section{Remarks and open problems}\label{sec:problems}

\begin{remark}[Triangular maps of the plane]\label{rem:R2}
Let $F(x,y)=(f(x),g(x,y))$ be a continuous triangular map of $\R^2$ and let $z\in\R^2$ have a bounded forward orbit with
closure $A$. Let $J$ and $L$ be the convex hulls of the coordinate projections of $A$ and $r_J$, $r_L$ the nearest-point
retractions of $\R$ onto them. Then $\tilde F(x,y)=(r_J(f(x)),r_L(g(x,y)))$ is a triangular map of $J\times L$ that agrees with
$F$ on $A$, so $\omega_{\tilde F}(z)=\omega_F(z)$. When $J$ and $L$ are nondegenerate, an affine change of coordinates
conjugates $\tilde F$ to a triangular map of $\I^2$; when one of them is a point, $\omega_F(z)$ lies in a line and the
one-dimensional results apply. Conversely, every triangular map of $\I^2$ extends to one of $\R^2$. Hence the results of
this paper hold for bounded orbits of triangular maps of $\R^2$, with $\I$ replaced by compact intervals.
\end{remark}

\begin{problem}\label{prob:W}
Prove or disprove Conjecture~\ref{conj} for tuples of sets in $\W$, that is, show whether condition (ii) alone implies
realizability for such tuples. Partial answers, and a description of the configurations that remain open, are given in
\cite{GuII}.
\end{problem}

\begin{problem}\label{prob:infinite}
When $\omega_f(x_0)$ is infinite, describe the possible multifunctions $x\mapsto\omega_F(z)\cap\I_x$ over $\omega_f(x_0)$.
Theorem~\ref{thm:A} suggests replacing the cyclic map by a map over the minimal set or $\omega$-limit set of the base.
\end{problem}

\begin{problem}\label{prob:hausdorff}
For triangular maps whose $\omega$-limit sets lie over a fixed periodic orbit of the base, describe when the family of
$\omega$-limit sets is closed in the Hausdorff metric. Already for one fibre it need not be \cite[Ch.~4]{Gu}, in contrast with
interval maps \cite{BBHS}.
\end{problem}

\begin{problem}\label{prob:unbounded}
Characterize the nonempty unbounded $\omega$-limit sets of triangular maps of $\R^2$.
\end{problem}


\end{document}